\input mssymb

\documentstyle[11pt]{article}

\title{$Aut(M)$ has a large dense
free subgroup
 for saturated $M$}

\author{Garvin Melles\thanks{Would like to thank Ehud Hrushovski
for supporting him with funds from NSF Grant DMS 8959511 and Wilfrid
Hodges for helping with the conjecture's history} \\Hebrew University of Jerusalem 
\and Saharon Shelah\thanks{partially supported by the U.S.-Israel Binational
Science Foundation. Publ. 452}\\Hebrew
University of Jerusalem\\Rutgers University}

\newtheorem{theorem}{Theorem}
\newtheorem{defi}[theorem]{Definition}
\newtheorem{lemma}[theorem]{Lemma}

\newcommand{\proof}{{\sc proof} \hspace{0.1in}}
\newcommand{\iopp}{{\bigcup}{\!\!\!\!{|}}}

\newcommand{\indep}[1]{\mathop{\iopp}\limits_{\textstyle{#1}}\ }
\newcommand{\sub}{\subseteq}

\begin{document}
\mathsurround=.1cm
\maketitle

\begin{abstract}
We prove that for a stable theory $T,$ if $M$ is a saturated model of
$T$ of cardinality $\lambda$ where $\lambda > \big|T\big|,$ then
$Aut(M)$ has a dense free subgroup on $2^{\lambda}$ generators. 
This affirms a conjecture of Hodges.
\end{abstract}

\section{Introduction}

A subgroup $G$ of the automorphism group of a model $M$ is said to be
dense if every finite restriction of an automorphism of $M$ can be extended to an
automorphism in $G.$ In this paper we present Shelah's proof of a conjecture of Hodges
that for a cardinal $\lambda$ with $\lambda > 
\big|T\big|,$ if $M$ is a saturated
model of $T$ of size $\lambda$ then the automorphism group of $M,$
$Aut(M),$ has a dense
free subgroup of cardinality $2^{\lambda}.$ 
Wilfrid Hodges had noted that the theorem was true for  $\lambda$
such that $\lambda \geq \big|T\big|$ and $\lambda^{<\lambda}=\lambda.$
Peter Neumann then pointed out to him that de Bruijn had shown that independently of
any set theoretic assumptions on $\lambda,$ $Sym(\lambda),$ the group
of permutations of $\lambda,$ has  
a free subgroup on $2^{\lambda}$ generators. On checking the proof, 
Hodges found one could also make the subgroup dense, so the natural
conjecture was that for any cardinal $\lambda>|T|$ if there is a
saturated model $M$ of cardinality $\lambda,$ then $Aut(M)$ has a dense
free subgroup on $2^{\lambda}$ generators. As Shelah likes
questions, Hodges asked him about the conjecture when Shelah went to
work with Hodges in the Summer of 1989.  The proof presented
in this paper is  simpler than the one from 1989, thought of by Shelah while he was helping
Melles with his earlier proof. While the proof is not complicated, Melles
thinks it is a nice application of the non forking relation for stable theories. By the following theorem of Shelah, the
only open case was for $T$ stable, although for completeness, we also
include here a proof for the case that 
$\big|M\big| = \lambda = \lambda^{<\lambda}.$

\begin{theorem}
$T$ has a saturated model in $\lambda$ iff one of the following hold
\begin{enumerate}
\item $\lambda = \lambda^{<\lambda}+\big|D(T)\big|$
\item $T$ is stable in $\lambda$
\end{enumerate}
\end{theorem}
\proof [Sh c] VIII 4.7.

\vspace{.1in}

\begin{defi}
A subgroup $G$ of $Aut(M)$ is $<\lambda$ dense if every elementary permutation of
a subset $A$ of $M$ such that $\big|A\big|<\lambda$ has an extension
in $G.$
\end{defi}

Melles asked Shelah about the natural stregthening of Hodges
question; Can one find a subroup $G$ of $Aut(M)$ for $M$ a saturated model of
cardinality $\lambda,$ such that $G$ is  $<\lambda$ dense and of
cardinality $2^{\lambda}?$
Shelah quickly found a proof affirming this stronger conjecture.

\noindent Throughout this paper we work in ${\frak C}^{eq}$ and follow the
notation from [Sh c]. See there for the definitions of any
notions left undefined here. We denote the identity map by $id.$

\section{$\big|M\big| = \lambda = \lambda^{<\lambda}$}

\begin{lemma}
Let $I$ be an infinite order and let $\langle a_i\mid i\in I\rangle$
be a sequence  indiscernible over 
$A,$ $f$ an elementary map with domain $A\cup \bigcup\Big\{a_i\mid
i\in I\Big\},$ and $B$ a set. Then there is a sequence $\langle
c_i\mid i\in I\rangle$ which realizes 
$tp\big(\langle f(a_i)\mid i\in I\rangle /f(A)\big)$ such that $\langle
c_i\mid i\in I\rangle$ is indiscernible over $B\cup f(A).$ 
\end{lemma}
\proof By compactness and Ramsey's theorem.

\begin{lemma}\label{free}
Let $\tau  = f_{n}^{\epsilon_n}\ldots f_{0}^{\epsilon_0}$ be a
term (intended to represent a composition of functions with $f^{1}$
meaning $f$ and $f^{-1}$ meaning the inverse of $f$) such that 
$\epsilon_i\in\{-1,1\}$ and $f_i = f_{i+1}\ \Rightarrow\
\epsilon_i = \epsilon_{i+1}.$ Let $M, N$ be models such that $M$ is saturated of
cardinality $\lambda,$ and $M\prec N$ with $N$ being $\lambda^+$
saturated and $\lambda^+$ homogenous. Let $\Big\{f_{\nu_0}\ldots
f_{\nu_n} \Big\}$ be a finite set of automorphisms of $M$ with
$f_{\nu_i}= f_{\nu_{i+1}}$ iff $f_i = f_{i+1}$ in $\tau.$ Then there
are automorphisms of $N,$ $\Big\{F_{\nu_0}\ldots
F_{\nu_n} \Big\}$ such that each $F_{\nu_i}$ extends $f_{\nu_i}$ and $F_{\nu_n}^{\epsilon_n} \ldots
F_{\nu_0}^{\epsilon_0} \neq id_N.$   
\end{lemma} 
\proof If $\epsilon_0=1,$ let $A_0=\big\{a_i^0\mid i<\omega\big\}\sub
N$ be an infinite indiscernible sequence over $M.$ Let $F$ be an
extension of $f_{\nu_0}$ to an automorphism of $N$ and let $A_1\sub N$
realize $tp(F[A_0]/F[M])$ such that $A_1$ is indiscernible over
$A_0\cup M.$ Let $F^0_{\nu_0}$ be the elementary map extending
$f_{\nu_0}$ such that $A_0$ is sent to $A_1.$ If $\epsilon_0=-1,$ then
let $A_0=\big\{a_i^0\mid i<\omega\big\}\sub 
N$ be an infinite indiscernible sequence over $M.$ Let $F$ be an
extension of $(f_{\nu_0})^{-1}$ to an automorphism of $N$ and let $A_1\sub N$
realize $tp(F[A_0]/F[M])$ such that $A_1$ is indiscernible over
$A_0\cup M.$ Let $(F^0_{\nu_0})^{-1}$ be the elementary map extending
$(f_{\nu_0})^{-1}$ such that $(F^0_{\nu_0})^{-1}$ sends $A_0$ to $A_1.$
Now by induction on $0<i\leq n$ we define infinite sequences $A_{i+1}$
indiscernible over $M\cup\bigcup\limits_{j\leq i}A_j$ and elementary
maps $F_{\nu_i}^i$ such that 
\begin{enumerate}
\item $F_{\nu_i}^i$ extends $f_{\nu_i}$
\item If $j<i$ and $\nu_j=\nu_i,$ then $F_{\nu_j}^j\sub F_{\nu_i}^i$
\item $F_{\nu_i}^{\epsilon_i} \ldots
F_{\nu_0}^{\epsilon_0}(A_0)=A_{i+1}$
\end{enumerate} 

\noindent Now suppose $0<i\leq n$ and $F^j_{\nu_j}$ have been defined
for all $j<i.$ Suppose $\epsilon_i=1$ and there is a $j<i$ such that
$\nu_j=\nu_i.$ Let $j^*$ be the largest such $j.$ Let $F$ be an
extension of $F^{j^*}_{\nu_{j^*}}$ to an automorphism of $N.$ By the
construction $A_i$ is indiscernible over the domain of
$F^{j^*}_{\nu_{j^*}}.$ So we can find $A_{i+1}$ realizing
$tp(F[A_i]/dom\,F^{j^*}_{\nu_{j^*}})$ such that $A_{i+1}$ is
indiscernible over $M\cup\bigcup\limits_{j\leq i}A_j.$ Let
$F_{\nu_i}^i$ be the
elementary map extending $F^{j^*}_{\nu_{j^*}}$ taking $A_i$ to $A_{i+1}.$ 
If $\epsilon_i=-1$ or if $j^*$ does not exist, the induction step is similar.
Now let $F_{\nu_j}$ be an automorphism of $N$ extending $F_{\nu_i}^i$
where $i$ is the largest index such that $\nu_i=\nu_j.$ 
$F_{\nu_n}^{\epsilon_n} \ldots
F_{\nu_0}^{\epsilon_0} \neq id_N$ since $F_{\nu_n}^{\epsilon_n} \ldots
F_{\nu_0}^{\epsilon_0}(A_0)=A_{n+1}$  and $A_0\cap A_{n+1}=\emptyset$ since $A_{n+1}$
is indiscernible over  $A_0\cup M.$

\begin{theorem}
Let $T$ be a complete theory, $M$ a saturated model of $T$ of cardinality
$\lambda$ with $\big|T\big| \leq \lambda = \lambda^{<\lambda}.$ Then
$Aut(M)$ has a dense free subgroup on $2^{\lambda}$ generators. 
\end{theorem}
\proof Let $TR =\, ^{<\lambda}\lambda.$ For $\alpha < \lambda,$ let
$TR_{\alpha} =\,^{<\alpha}\lambda$ and let $\bar{0}_{\alpha}$ be the
function with domain $\alpha$ and range $\big\{0\big\}.$ We define by
induction on $\alpha < \lambda$ a model $M_{\lambda}$ of $T$ and
$f_{\eta}\in Aut(M_{\alpha})$ for $\eta\in TR_{\alpha}$ such that
\begin{enumerate}
\item $M_{\alpha}\models T$
\item $\big|M_{\alpha}\big| =  \lambda$
\item $\langle M_{\alpha}\mid \alpha < \lambda\rangle$ is increasing
continuous
\item If $\alpha$ is a successor, then $M_{\alpha}$ is saturated
\item $\nu\lhd\eta\rightarrow f_{\nu}\subseteq f_{\eta}$
\item If $\alpha=\beta+1,$  then $\langle f_{\eta}\mid \eta\in TR_{\alpha}\backslash
\Big\{\bar{0}_{\beta}\frown i\mid  i<\lambda\Big\}\rangle$  is free
\end{enumerate}
For $\alpha=\beta +1$ we let $\langle f_{\bar{0}_{\beta}\frown i}\mid
i < \lambda\rangle$ be a sequence of automorphisms of $M_{\alpha}$
such that each finite partial automorphism of $M_{\beta}$ has an
extension in $$\Big\{ f_{\bar{0}_{\beta}\frown i}\mid
i < \lambda\Big\}.$$
If we succeed in doing the induction then for $\eta\in\ 
^{\lambda}\lambda$ if
$$f_{\eta}=\bigcup\Big\{f_{\nu}\mid \nu = \eta\restriction  \alpha,
\ \alpha < \lambda\Big\}$$
then the $f_{\eta}$  and the $\bigcup\limits_{\alpha <
\lambda}M_{\alpha}$ are as required by the theorem. The only problem
in the induction is at successor steps, so let $\alpha = \beta + 1.$
Let
$$\Gamma = \Big\{\tau\mid \tau\hbox{  is a term of the form
}f_{\nu_n}^{\epsilon_n}\ldots f_{\nu_0}^{\epsilon_0} \Big\}$$
such that
\begin{enumerate}
\item $\forall i  <  n+1,\  \epsilon_i\in \{-1,1\}$
\item $\forall i  <  n+1,\   \nu_i\in TR_{\alpha}\backslash
\Big\{\bar{0}_{\beta}\frown i\mid  i<\lambda\Big\}$
\item $\forall i  <  n+1,\  \nu_i = \nu_{i+1}\ \Rightarrow\ \epsilon_i =
\epsilon_{i+1}$ 
\end{enumerate}

\noindent Let $\langle \tau_i\mid i < \lambda\rangle$ be a well ordering of
$\Gamma.$ Let $N$ be a $\lambda^+$ saturated, $\lambda^+$ homogenous
model of $T$ containing $M_{\beta}.$ By induction on $\gamma<\lambda$
we define elementary submodels $M_{\beta,\gamma}$ of $N$ and for every
$\nu\in TR_{\alpha}\backslash
\Big\{\bar{0}_{\beta}\frown i\mid  i<\lambda\Big\}, \ f_{\nu,\gamma}$
such that
\begin{enumerate}
\item $M_{\beta,0}=M_{\beta}$
\item $f_{\nu,0}= f_{\nu\restriction\alpha}$
\item If $\gamma=\zeta+1,\  M_{\beta,\gamma}$ is saturated of cardinality $\lambda$
\item $\zeta<\gamma\ \Rightarrow\ f_{\nu,\zeta}\subseteq
f_{\nu,\gamma}$
\item If $\gamma=\zeta+1$ and
$\tau_{\zeta}=f_{\nu_n}^{\epsilon_n}\ldots f_{\nu_0}^{\epsilon_0}$
then $f_{\nu_n,\gamma}^{\epsilon_n}\ldots f_{\nu_0,\gamma}^{\epsilon_0}\neq
id_{M_{\beta,\gamma}}$
\end{enumerate} 

\noindent If we succeed in the induction then we can let
$M_{\alpha}=\bigcup\limits_{\alpha<\lambda}M_{\beta,\gamma}$ and
$f_{\nu}=\bigcup\limits_{\alpha<\lambda}f_{\nu,\gamma}$. 
The only non-trivial step in the induction is for successor steps, so
let $\gamma = \zeta +1.$ By lemma \ref{free} we can find automorphisms
$F_{\nu_0,\gamma},\ldots,F_{\nu_n,\gamma}$ of $N$ extending
$f_{\nu_0,\zeta},\ldots,f_{\nu_n,\zeta}$ such that if
$\tau_{\zeta} = f_{\nu_n}^{\epsilon_n}\ldots f_{\nu_0}^{\epsilon_0}$
then
$$F_{\nu_n,\gamma}^{\epsilon_n}\ldots
F_{\nu_0,\gamma}^{\epsilon_0}\neq id_N$$
For $\nu \in TR_{\alpha}\backslash
\Big\{\bar{0}_{\beta}\frown i\mid  i<\lambda\Big\},$ but not in
$\{\nu_0,\ldots,\nu_n\},$ let $F_{\nu,\gamma}$ be an arbritary
extension of $f_{\nu,\zeta}$ to $N.$ Let $M_{\beta,\gamma}\prec N$ be
a saturated model of size $\lambda$ such that 
\begin{enumerate}
\item $M_{\beta,\zeta}\prec M_{\beta,\gamma}$
\item $M_{\beta,\gamma}$ contains witnesses to $F_{\nu_n,\gamma}^{\epsilon_n}\ldots
F_{\nu_0,\gamma}^{\epsilon_0}\neq id_N$
\item $M_{\beta,\gamma}$ is closed under the $F_{\nu,\gamma}$
\end{enumerate}
For each $\nu\in TR_{\alpha}\backslash
\Big\{\bar{0}_{\beta}\frown i\mid  i<\lambda\Big\},$ let
$f_{\nu,\gamma}$ be $F_{\nu,\gamma}\restriction M_{\beta,\gamma}.$

\section{$\big| M\big| = \lambda < \lambda^{<\lambda}$}

\noindent Throughout this section, by  theorem 1 mentioned in the
introduction, we can assume that $T$ is stable. Although the proofs in this section
are  simple, there is an hidden element of complexity covered over by theorem 1.   
\begin{theorem}\label{chains}
Let  $\langle M_i\mid i<\delta\rangle$ is an increasing  
elementary chain
of models of $T$ that are $\lambda$ saturated  with $cf
\,\delta\geq \kappa_r(T).$ Then 
$\bigcup\limits_{i<\delta}M_i$ is a $\lambda$ saturated model of $T.$
\end{theorem}
\proof [Sh c] III 3.11

\begin{lemma}\label{indep}
Let $\Big\{C_i\mid i\in I\Big\}$ be independent over $A$ and let
$\Big\{D_i\mid i\in I\Big\}$ be independent over $B.$ Suppose that for
each $i\in I$, $tp(C_i/A)$ is stationary. Let $f$ be an elementary map
from $A$ onto $B,$ and let for each $i\in I,$ $f_i$ be an elementary
map extending $f$ which sends $C_i$ onto $D_i.$ Then 
$$\bigcup\limits_{i\in I}f_i$$
is an elementary map from $\bigcup\limits_{i\in I}C_i$ onto $\bigcup\limits_{i\in I}D_i.$
\end{lemma}
\proof Left to the reader.

\begin{theorem}
Let $T$ be a complete stable theory and let $M$ be a saturated model
of $T$ of cardinality $\lambda >  \big|T\big|.$ Then 
\begin{enumerate}
\item $Aut(M)$ has a dense free subgroup $G$ of cardinality $2^{\lambda}$
\item In fact, if $\sigma\leq\lambda$ is regular, then there is a free
subgroup $G$ of $Aut(M)$ such that any partial
automorphism $f$ of $M$ with $\big|dom\,f\big|<\sigma$ can be
extended to an element of $G.$
\end{enumerate}
\end{theorem}
\proof Let $\sigma+\kappa_r(T)\leq \kappa=cf(\kappa)\leq\lambda.$ 
We define by induction on $i\leq\kappa$ an increasing continuous
elementary chain of models $M_i$ of $T,$ ordinals $\alpha_i$ of
cardinality $2^{\lambda}$ and families $\Big\{g_{\alpha}^i\mid
\alpha<\alpha_i\Big\}$ of automorphisms of $M_i$ and such that

\begin{enumerate}
\item $\alpha_0=2^{\lambda}$
\item $\big|M_i\big|= \lambda$
\item $\langle \alpha_j\mid j\leq i\rangle$ is increasing continuous
\item $\forall g\in Aut(M_i)\ \ \ \bigvee_{\alpha<\alpha_i}g\subseteq 
g_{\alpha}^{i+1}$
\item For a fixed $\alpha,$ the $g_{\alpha}^i$ form an elementary chain
\item $\langle g_{\alpha}^{i+1}\mid \alpha < \alpha_{i+1}\rangle$ \ is free
\item If $i=j+1,$ or $i=0$ then $M_i$ is saturated.
\end{enumerate}

\noindent If we succeed in doing the induction then by theorem \ref{chains} $M_{\kappa} =
\bigcup\limits_{i<\kappa}M_i$ is a saturated model of cardinality $\lambda.$
If we let for $\alpha<\bigcup\limits_{i<\kappa}\alpha_i,$
$$g_{\alpha}=\bigcup\Big\{g^j_{\alpha}\mid \alpha_j > \alpha,
j<\kappa\Big\}$$ then $\Big\{g_{\alpha}\mid
\alpha<\bigcup\limits_{i<\kappa}\alpha_i\Big\}$ is free by item $6.$ in the
construction and is dense (in the strong sense of $2.$ of the theorem) by
item $5.$ 

\vspace{.1in}

\noindent The only difficulty in the induction is for $i=j+1.$ 
Let $\Big\{p_{\zeta}\mid \zeta < \zeta^*\Big\}$ list
$S^1(acl\,\emptyset).$ (So $\zeta^*\leq\lambda$) Let
$\Big\{a^{\zeta}_{\gamma}\mid \zeta<\zeta^*, \gamma < \lambda 
\Big\}$ be a set of elements independent over $M_j$ such that
$tp(a^{\zeta}_{\gamma}/M_j)$ is a nonforking extension of $p_{\zeta}.$ 
For every $g\in Aut(M_j)$ let $f^{[g]}$ be
the permutation of $\zeta^*$ such that
$$f^{[g]}(\zeta)=\xi\ \Leftrightarrow\ g(p_{\zeta})= p_{\xi}$$
List $Aut(M_j)\backslash\Big\{g_{\alpha}^j\mid \alpha<\alpha_j\Big\}$ as
$$\langle g^j_{\alpha}\mid \alpha_j\leq \alpha < \alpha_i\rangle$$
Let
$$A_i= M_j\cup\Big\{a^{\zeta}_{\gamma}\mid \zeta < \zeta^*, \gamma<
\lambda \Big\}$$
and let
$$\Big\{h_{\alpha}^i\mid \alpha < \alpha_i\Big\}$$
be a set of free permutations of $Sym(\lambda).$ Define for $\alpha
<\alpha_i$ a permutation $g^{j,*}_{\alpha}$ of $A_i$ by letting
$g^{j,*}_{\alpha}\restriction M_j=g_{\alpha}^j$ and
$$g^{j,*}_{\alpha}(a^{\zeta}_{\gamma}) = a^{f^{[g^j_{\alpha}]}(\zeta)}_{h_{\alpha}^i(\gamma)}$$
By lemma \ref{indep} each $g^{j,*}_{\alpha}$ is an elementary map.  
The $\langle g^{j,*}_{\alpha}\mid \alpha < \alpha_i\rangle$ are free.
For suppose not. Then for some
$\big\{\alpha_1,\ldots,\alpha_n,\alpha_{n+1}\big\}\subseteq \alpha_i$
we would have
$$g^{j,*}_{\alpha_1}\ldots g^{j,*}_{\alpha_n}=
g^{j,*}_{\alpha_{n+1}}$$
so for every $\zeta<\zeta^*$
$$f[g^{j,*}_{\alpha_1}\ldots g^{j,*}_{\alpha_n}](\zeta)= f[g^{j,*}_{\alpha_{n+1}}](\zeta)$$
and for every $\gamma < \lambda,$
$$h_{\alpha_1}\ldots h_{\alpha_n}(\gamma)= h_{\alpha_{n+1}}(\gamma)$$
a contradiction to the freeness of the $h_{\alpha_i}.$
Let $M_i$ be a model such that
\begin{enumerate}
\item $A_i\subseteq M_i \prec {\frak C}^{eq}$
\item $\big(M_i,a\big)_{a\in A_i}$ is saturated
\item $\big|M_i\big|=\lambda$
\end{enumerate}
This is possible as the theory of $\big({\frak C}^{eq},a\big)_{a\in
A_i}$ is stable in $\lambda.$ As $\big(M_i,a\big)_{a\in A_i}$ is
saturated we can  for each $\alpha<\alpha_i,$
let $g^i_{\alpha}$ be an extension of $g^{j,*}_{\alpha}$ to an automorphism
of $M_i.$

\section{$<\lambda$ Denseness}

\begin{theorem}
Let $M$ be a saturated model of cardinality $\lambda>|T|.$ Then
$Aut(M)$ has a free $<\lambda$ dense free subgroup on $2^{\lambda}$
generators. 
\end{theorem}
\proof If $\lambda^{<\lambda}=\lambda,$ the proof given gives a
$<\lambda$ dense free subgroup. So we can assume that $T$ is stable in $\lambda.$
We work in ${\frak C}^{eq}.$ Let $\langle p_i\mid i<i^*\leq
\lambda\rangle$ list all types over $acl\,\emptyset.$ Let $\Big\{
a_{i,\zeta,\xi}\mid i<i^*,\ \zeta<\lambda,\ \xi<\lambda\Big\}$ be
independent over $\emptyset,$ with $a_{i,\zeta,\xi}$ realizing $p_i$
and 
$$\big(M,a_{i,\zeta,\xi}\big)_{(i,\zeta,\xi)\in
i^*\times\lambda\times\lambda}$$ 
is saturated. 
Let $\Big\{f_{\alpha}\mid \alpha<2^{\lambda}\Big\}$ be a free subgroup
of $Sym(\lambda).$ Let $\Big\{g_{\alpha}\mid \alpha<2^{\lambda}\Big\}$
be a list of permutations of subsets of $M$ of cardinality $<\lambda$
such that for every $\alpha<2^{\lambda},$ $acl\,\emptyset\subseteq
dom\,g_{\alpha}.$ 
Let $C_{\alpha}= dom\,g_{\alpha}\  (=ran\,g_{\alpha}).$ For some
subset $u_{\alpha}$ of $i^*\times\lambda\times\lambda$ such that
$\big|u_{\alpha}\big|\leq\big|C_{\alpha}\big|+\kappa_r(T),$
$$C_{\alpha}\,\indep{\big\{a_{i,\zeta,\xi}\mid (i,\zeta,\xi)\in
u_{\alpha}\big\}}\,\big\{a_{i,\zeta,\xi}\mid (i,\zeta,\xi)\in
i^*\times\lambda\times\lambda\big\}$$
We can find a $D_{\alpha}\supseteq C_{\alpha}$ and
$v_{\alpha}\supseteq u_{\alpha}$ such that
$\big|D_{\alpha}\big|=\big|C_{\alpha}\big|,\ \big|v_{\alpha}
\big|\leq \big|C_{\alpha}\big|+\kappa_r(T),$ and for some extension $g_{\alpha}'$ of
$g_{\alpha},$ $g_{\alpha}'$ is an automorphism of $D_{\alpha}$ with
$D_{\alpha}\supseteq \big\{a_{i,\zeta,\xi}\mid (i,\zeta,\xi)\in
v_{\alpha}\big\}$ and
$$D_{\alpha}\,\indep{\big\{a_{i,\zeta,\xi}\mid (i,\zeta,\xi)\in
v_{\alpha}\big\}}\,\big\{a_{i,\zeta,\xi}\mid (i,\zeta,\xi)\in
i^*\times\lambda\times\lambda\big\}$$
Since
$$\big\{a_{i,\zeta,\xi}\mid (i,\zeta,\xi)\in
i^*\times\lambda\times\lambda-v_{\alpha}\big\}\,\indep{\emptyset}\,
\big\{a_{i,\zeta,\xi}\mid
(i,\zeta,\xi)\in  
v_{\alpha}\big\}$$
we have
$$D_{\alpha}\,\indep{\emptyset}\,\big\{a_{i,\zeta,\xi}\mid (i,\zeta,\xi)\in
i^*\times\lambda\times\lambda-v_{\alpha}\big\}$$
For each $\alpha<2^{\lambda}$ let 
$$G_{\zeta}^{\alpha}=\big\{\zeta<\lambda\mid \forall\,\xi<\lambda\
\forall\,i<i^*\ a_{i,\zeta,\xi}\not\in D_{\alpha}\big\}$$
Let $h_{\alpha}$ be the map taking $a_{i,\zeta,\xi}$ to
$a_{i',\zeta,f_{\alpha}(\xi)}$ for $\zeta\in G^{\alpha}_{\zeta}$ if
$g_{\alpha}(p_i)=p_{i'}.$  Since 
$$D_{\alpha}\,\indep{\emptyset}\,\big\{a_{i,\zeta,\xi}\mid (i,\zeta,\xi)\in
i^*\times\lambda\times\lambda-v_{\alpha}\big\}$$
and $g_{\alpha}'$ and $h_{\alpha}$ agree on $acl\,\emptyset,$
$g_{\alpha}'\cup h_{\alpha}$ is an elementary map. Let $g_{\alpha}''$
be an extension of $g_{\alpha}'\cup h_{\alpha}$ to an automorphism of
$M.$ (This is possible as
$$\big(M,c\big)_{c\in D_{\alpha}\cup dom\,h_{\alpha}\cup
ran\,h_{\alpha}}$$
is saturated.) If $\big\{\alpha_0,\ldots,\alpha_n\big\}\subseteq
2^{\lambda}$ then 
$$G_{\zeta}^{\alpha_0}\cap\ldots\cap G_{\zeta}^{\alpha_n}\neq \emptyset$$
so the $g_{\alpha}''$ are free, and by construction $g_{\alpha}''$
extends $g_{\alpha}$ so the $g_{\alpha}''$ are $<\lambda$ dense.

\pagebreak

\begin{center}
REFERENCES
\end{center}

\begin{enumerate}

\item N.G. de Bruijn, {\em Embedding Theorems for Infinite Groups},
Nederl. Akad. Wetensch. Indag. Math. 19 (1957), p. 560-569.

\item Richard Kaye, {\em The automorphism group of a countably
recursive saturated model, } Proc. London Math Soc. (3) 65 (1992),
p. 225-244.

\item D. Macpherson, {\em Groups of automorphisms of a $\aleph_0$
categorical structure,} Quarterly J. Math. Oxford 37 (1986) p.
449-465.

\item [Sh c] Saharon Shelah, {\em Classification theory and the
number of isomorphic models, revised}, North Holland Publ. Co. Amsterdam,
Studies in Logic and the foundations of Math.,  vol  92, 1990.

\end{enumerate}

\end{document}